# ON FINITE TYPE 3-MANIFOLD INVARIANTS II

STAVROS GAROUFALIDIS AND JEROME LEVINE

Fax number for J. Levine: (617) 736 3085

Email: stavros@math.mit.edu

Email: levine@max.math.brandeis.edu

Abstract. The present paper is a continuation of [Oh] and [Ga] devoted to the study of finite type invariants of oriented integral homology 3-spheres ($\mathbb{Z}HS$ for short). The paper consists of two parts. In the first part we classify pure braids and string links modulo the relation of surgical equivalence. We prove that the group of surgical equivalence classes of pure braids is isomorphic to the corresponding group of string links (theorem 2). We also give two alternative descriptions of the above mentioned group $P^{\mathbb{SE}}(n)$ of surgical equivalence classes of $n$ strand pure braids: one as a semidirect product of $P^{\mathbb{SE}}(n-1)$ together with an explicit quotient of the free group. And another description (theorem 3) as a group of automorphisms of a nilpotent quotient of a free group. In the second part we study the finite type invariants of $\mathbb{Z}HS$, originally introduced by Ohtsuki [Oh] and partially answer questions 1 and 2 from [Ga]. We reprove Ohtsuki's fundamental result which states that the space of type $m$ invariants of $\mathbb{Z}HS$ is finite dimensional for every $m$. Our proof allows us to show (corollary 3.5) that the graded space of degree $m$ invariants of $\mathbb{Z}HS$ is zero dimensional unless $m$ is divisible by 3. This partially answers question 1 of [Ga]. Furthermore, we study a map from knots (in $S^3$) to $\mathbb{Z}HS$, and show that type $5m+1$ invariants of $\mathbb{Z}HS$ map to type $4m$ invariants of knots, thus making progress towards question 2 of [Ga].

## Contents

The authors were partially supported by NSF grants DMS-95-05105 and DMS-93-03489.

# 1. Introduction

**1.1. History.** The present paper is a continuation of [Oh] and [Ga], devoted to the study of finite type invariants of *oriented* integral homology 3-spheres ($\mathbb{Z}HS$ for short). Finite type invariants of $\mathbb{Z}HS$ were originally introduced by Ohtsuki [Oh] in his seminal paper (for a precise definition, see definition 3.2). There are at least two sources of motivation/analogy/inspiration: the finite type knot invariants (see [B-N], [BL], [Va]) and Chern-Simons theory in 3 dimensions (see [Wi2], [Rz1], [Rz2]).

Recall that a type $m$ (otherwise called Vassiliev) invariant of knots in $S^3$ (for a precise definition see the references above) satisfies a difference formula with respect to cutting along $m+1$ spheres, twisting and gluing back. Vassiliev invariants of knots in $S^3$ form a unifying way of thinking about the Alexander, Jones, HOMFLY, Kauffman (and others) polynomials of knots.

Similarly, a type $m$ invariant of $\mathbb{Z}HS$ satisfies a "difference formula" with respect to cutting $m+1$ solid torii from a $\mathbb{Z}HS$, twisting them, and gluing them back. For a precise statement, see definition 3.2. It is hoped that finite type invariants of $\mathbb{Z}HS$ will provide a unifying way of dealing with 3-manifold invariants and with Chern-Simons theory in a way analogous to the case of Vassiliev invariants.

In either of the above sources of inspiration, we can think of type $m$ invariants of $\mathbb{Z}HS$ as "polynomials of degree $m$" on the infinite dimensional vector space $\mathcal{M}$. Here and in the rest of this paper $\mathcal{M}$ is the vector space (over $\mathbb{Q}$) with basis the set of $\mathbb{Z}HS$. It is an important question to ask whether finite type invariants of $\mathbb{Z}HS$ separate points in $\mathcal{M}$.

**1.2. Statement of results; plan of the proof.** This paper consists of two parts. In the first part, which consists of section 2 we classify pure braids and string links modulo the relation of surgical equivalence. In section 2.1 we recall some definitions of string links and pure braids. In section 2.2 we study a subgroup $A(F^4(n))$ of the automorphism group of the nilpotent quotient $F^4(n)$ of a free group and show in theorem 1 how to express $A(F^4(n))$ as a semidirect product of $A(F^4(n-1))$ together with a nilpotent quotient of a free group. In section 2.3 we prove that the group of surgical equivalence classes of pure braids is isomorphic to the corresponding group of string links (theorem 2). We also give two alternative descriptions of the above mentioned group $P^{\mathbb{SE}}(n)$ of surgical equivalence classes of $n$ strand pure braids: one as a semidirect product of $P^{\mathbb{SE}}(n-1)$ together with an explicit quotient of the free group and another description (theorem 3) as a group of automorphisms of a nilpotent quotient of a free group.

In the second part which consists of section 3 we study the finite type invariants of $\mathbb{Z}HS$, originally introduced by Ohtsuki [Oh] and partially answer questions 1 and 2 from [Ga]. In section 3.2 we reprove Ohtsuki's fundamental result which states that the space of type $m$ invariants of $\mathbb{Z}HS$ is finite dimensional for every $m$. Our proof allows us to show a vanishing theorem (corollary 3.5), namely that the graded space of degree $m$ invariants of $\mathbb{Z}HS$ is zero dimensional unless $m$ is divisible by 3. This partially answers question 1 of [Ga]. In section 3.4 we study the map from knots (in $S^3$) to $\mathbb{Z}HS$ defined by mapping a knot $K$ in $S^3$ to the $\mathbb{Z}HS$ $(S^3)^{K,+1}$ obtained by $+1$ surgery on $K$ (see [Ro], and section 3.1). We show (proposition 3.9) that type $5m+1$ invariants of $\mathbb{Z}HS$ map to type $4m$ invariants of knots, thus making progress towards answering question 2 of [Ga]. We end up in section 4 with a philosophical comment about the appearance of trivalent graphs in this paper.

After the present work was completed we received a copy of [GL] in which they prove a very special case of corollary 3.5 as well as a weak form of proposition 3.9.

**1.3. Acknowledgment.** We wish to thank D. Bar-Natan for for many useful conversations. Especially we wish to thank T. Ohtsuki for enlightening electronic email conversations and helpful comments, and the `Internet` for providing the support for the relevant communications.

## 2. Surgical equivalence of string links and pure braids

In this section we give two different classification theorems of the group of surgical equivalence classes of sting links and pure braids. As in the case of string link homotopy (see [HL]) these groups are isomorphic. These classification theorems are entirely analogous to the classical theorems of Artin [Ar] on isotopy classification of braids and the more recent results of Habegger-Lin [HL] on homotopy classification of string links. The first theorem is a recursive description of the $n$ strand pure braid group in terms of the $n-1$ strand one and a certain nilpotent quotient of a free group. The second theorem gives an isomorphism with a certain group of automorphisms of another nilpotent quotient of a free group.

The proof of the first theorem is similar to that of [HL] while the second is rather different (*warning*: lemma 1.9 of [HL] is not correct as stated and the proof of theorem 1.7 of that paper seems to require some modification).

**2.1. Preliminaries on surgical equivalence.** We begin by recalling some definitions.

**Definition 2.1.** (1) Let $n$ be a positive integer, $I$ be the unit interval $[0,1]$ and $D^2$ be the unit disc in the plane. An *$n$-string link* is a disjoint union of $n$ (smooth) arcs $\sigma_1, \cdots, \sigma_n$ in the solid cylinder $I \times D^2$ so that the boundary of the $i$-th arc is $\partial I \times p_i$, where $\{p_i\}$ are $n$ distinct points in $D^2$.

(2) An *$n$-strand pure braid*[1] is an $n$-string link with the extra property that the tangent vector at any point of any $\sigma_i$ is never horizontal.See figure 1.

(3) A *surgery* on a string link $\sigma$ produces another string link $\sigma'$ as follows. Let $\gamma$ be an unknotted closed curve in $I \times D^2 - \sigma$. If we remove a tubular neighborhood of $\gamma$ and sew it back in so that the new meridian is identified with a former longitude which links $\gamma$ once, then $I \times D^2$ is converted into a new manifold $\Delta$ which is diffeomorphic to $I \times D^2$ again. We define $(I \times D^2, \sigma') = (\Delta, \sigma)$. A more concrete description of $\sigma'$ is given by removing a tubular neighborhood of a disk bounded by $\gamma$ and then reinserting it with a single clockwise or counterclockwise twist. If $\sigma$ is a braid then so is $\sigma'$. We say two string links or pure braids are *surgically equivalent* if one can be obtained from the other by a sequence of surgeries. It is clear that homotopic string links or braids (i.e., string links for which there is a homotopy $h_t$ such that distinct components remain disjoint for all $t \in [0,1]$. See also [Mi1], [Mi2], [HL] and figure 4.) are surgically equivalent. If $\sigma_1, \sigma_2$ are two string links, then the product string link $\sigma_1\sigma_2$ is obtained by stacking $\sigma_1$ on top of $\sigma_2$. If $\overline{\sigma}$ is the reflection of $\sigma$ about $\frac{1}{2} \times D^2$, then $\sigma\overline{\sigma}$ and $\overline{\sigma}\sigma$ are homotopic and, therefore, surgically equivalent to the trivial string link (see [HL]). Thus we obtain groups $P^{\mathbb{SE}}(n)$ (respectively,

[1] we apologize for using the words strand (for pure braids) and string (for string links) to mean components of links of braids

$SL^{\mathbb{SE}}(n)$) of surgical equivalence classes of $n$-strand pure braids (repectively, $n$-string links). There is an obvious homomorphism $P^{\mathbb{SE}}(n) \to SL^{\mathbb{SE}}(n)$.

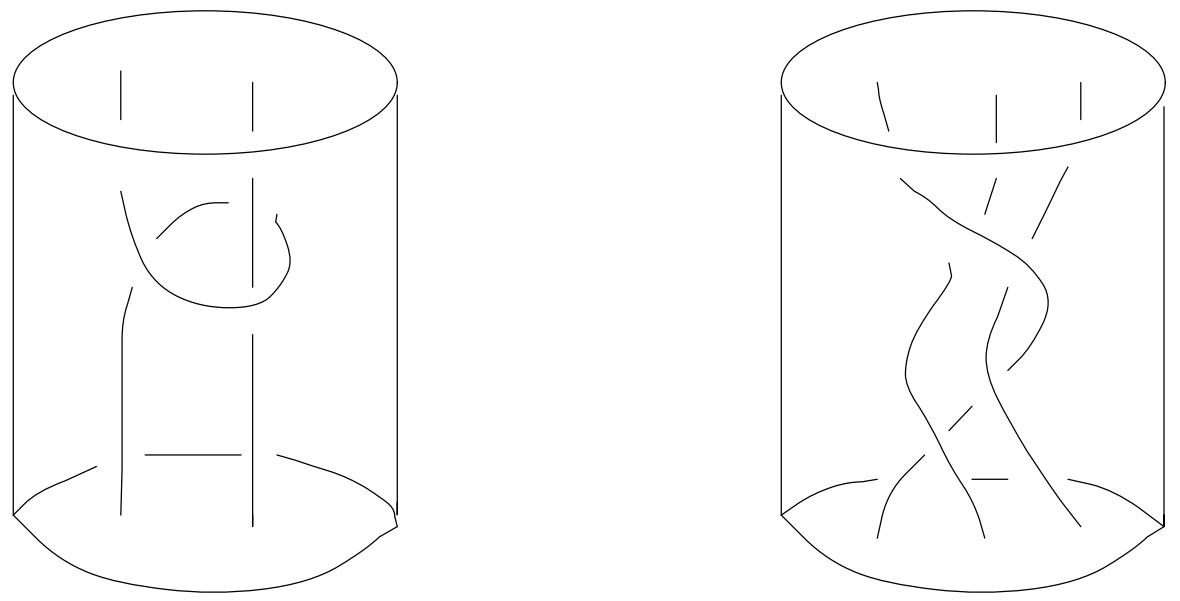

**Figure 1.** A string link of 2 components (on the left) and a braid of 3 strands on the right

Let $F(m)$ denote the free group on generators $\{x_1, \cdots, x_m\}$. Recall the homomorphism $\tau : F(n-1) \to P^{\mathbb{I}}(n)$ defined by inserting the $n^{th}$ strand into a trivial $(n-1)$-strand braid. Here $P^{\mathbb{I}}(n)$ stands for the group of isotopy classes of pure braids in $n$ strands. For any group $G$ let $G_q$ denote the subgroup generated by all commutators of order $q$, with the convention that $G_1$ is equal to $G$ and $G_2 = [G, G]$. Let $G^q$ denote the quotient group $G/G_q$.

**Claim 2.2.** *For any $w \in F(n-1)$ the surgical equivalence class of $\tau(w)$ depends only on the image of $w$ in $F^3(n-1)$.*

This follows directly from [Le1].

Thus we have a homomorphism $\tau : F^3(n-1) \to P^{\mathbb{SE}}(n)$. We also recall the obvious split epimorphisms $P^{\mathbb{SE}}(n) \to P^{\mathbb{SE}}(n-1)$ and $SL^{\mathbb{SE}}(n) \to SL^{\mathbb{SE}}(n-1)$ defined by deleting the $n^{th}$ strand. Let $\tau' : F^3(n-1) \to SL^{\mathbb{SE}}(n)$ be the composition $F^3(n-1) \xrightarrow{\tau} P^{\mathbb{SE}}(n) \longrightarrow SL^{\mathbb{SE}}(n)$.

### 2.2. The study of $A(F^4(n))$.

In this section we study a subgroup $A(F^4(n))$ of the group of automorphisms of $F^4(n)$ and prove theorem 1. Theorem 1 is the key ingredient in determining the structure of the groups $P^{\mathbb{SE}}(n)$ and $SL^{\mathbb{SE}}(n)$ in the next section.

First consider all automorphisms $\alpha$ which send each $x_i$ to a conjugate of itself. In fact, for any sequence of elements $g_1, \cdots, g_n$ of $F(n)$, there is an automorphism $\alpha$ such that $\alpha(x_i) = g_i x_i g_i^{-1}$. These equations only define endomorphisms of $F(n)$, in general, but always define automorphisms of $F^q(n)$ for any $q \geq 2$. It is easy to see that $\alpha$ depends only on the class of the $\{g_i\}$ in $F^{q-1}(n)$. If we demand that $g_i$ is *$i$-reduced*, i.e. the exponent sum of $x_i$ in $g_i$ is zero, then the $g_i \in F^{q-1}(n)$ are determined by $\alpha$. We define $A(F^q(n))$ to be the subgroup of all such $\alpha$ satisfying the additional property:

$$\alpha(x_1 \dots x_n) = x_1 \dots x_n \tag{1}$$

The following lemma will be useful in our study of the group $A(F^4(n))$.

**Lemma 2.3.** (a) *If* $\lambda \in F^4_2(n) \cap \langle x_n \rangle$*, then there exist unique elements* $\lambda_i \in F^3(n) \cap \langle x_n \rangle$ *such that*

$$[x_1, \lambda_1] \cdots [x_n, \lambda_n] = \lambda \tag{2}$$

*where* $\langle x_n \rangle$ *denotes the normal closure of* $x_n$ *in* $F^4_2(n)$*.*

(b) *Suppose* $\lambda_1, \cdots, \lambda_n \in F^3(n-1)$ *satisfy*

$$\lambda_1 x_1 \lambda_1^{-1} \cdots \lambda_{n-1} x_{n-1} \lambda_{n-1}^{-1} = x_1 \cdots x_{n-1} \quad \text{in } F^4(n-1) \tag{3}$$

*Then there exists unique* $\tilde{\lambda}_1, \cdots, \tilde{\lambda}_n \in F^3(n)$ *such that* $\tilde{\lambda}_i \equiv \lambda_i \mod \langle x_n \rangle$, $\tilde{\lambda}_n$ *is n-reduced and*

$$\tilde{\lambda}_1 x_1 \tilde{\lambda}_1^{-1} \cdots \tilde{\lambda}_n x_n \tilde{\lambda}_n^{-1} = x_1 \cdots x_n \quad \text{in } F^4(n) \tag{4}$$

*Remark* 2.4. This lemma is *false* if we replace $F^4(n)$ and $F^3(n)$ by $F^q(n)$ and $F^{q-1}(n)$ if $q > 4$. Also, note that, in (b), $\tilde{\lambda}_i$ is $i$-reduced if and only if $\lambda_i$ is $i$-reduced.

*Proof.* To prove part $(a)$ recall the notion of a Hall basis (see [MKS]). In particular, any element of $F^4(n)$ can be uniquely written in the following form:

$$x_1^{e_1} \cdots x_n^{e_n} \prod_{i<j} [x_i, x_j]^{e_{ij}} \prod_{\substack{i<j \\ k \leq j}} [x_k, [x_i, x_j]]^{e_{ijk}} \tag{5}$$

Thus we may uniquely write:

$$\lambda = \prod_{i,n} [x_i, x_n]^{e_i} \prod_{\substack{i<n \\ j \leq n}} [x_j, [x_i, x_n]]^{e_{ij}} \tag{6}$$

Therefore, equation (2) has a solution as follows:

$$\lambda_j = x_n^{e_j} \prod_{i<n} [x_i, x_n]^{e_{ij}} \tag{7}$$

We have set $e_n = 0$. Since any allowable solution $\lambda_j$ has such a representation, uniqueness follows also.

In order to show part $(b)$, we observe that any acceptable solution $\{\tilde{\lambda}_i\}$ of equation (4) can be written in the form $\tilde{\lambda}_i = \lambda_i \alpha_i x_n^{e_i}$ for some $\alpha_i \in F^3_2(n)$, with $e_n = 0$. Substituting this into equation (4) with a little bit of commutator manipulation and using the fact that $F^4_3(n)$ is the center of $F^4(n)$, we obtain the following formula in $F^4(n)$:

$$\left(\prod_{i=1}^{n-1} [\lambda_i, [x_n^{e_i}, x_i]]\right)\left(\prod_{i=1}^{n} [\alpha_i, x_i]\right) \prod_{i=1}^{n} ([x_n^{e_i}, x_i] \lambda_i x_i \lambda_i^{-1}) = \prod_{i=1}^{n} x_i \tag{8}$$

Reducing to $F^3(n)$ we obtain the much simpler formula:

$$(\prod_{i=1}^{n-1}[x_n, x_i]^{e_i})\prod_{i=1}^{n}\lambda_i x_i \lambda_i^{-1} = \prod_{i=1}^{n} x_i \tag{9}$$

Substituting from equation (3) this becomes:

$$\prod_{i=1}^{n-1}[x_n, x_i]^{e_i} = [x_n, \lambda_n] \tag{10}$$

and so $\lambda_n$ determines the $\{e_i\}$. Now that the $\{e_i\}$ are known, equation (8) can be written as $\prod_{i=1}^{n}[\alpha_i, x_i] = \tau$, where $\tau$ is a specified element of $F^4(n) \cap \langle x_n \rangle$. It follows from (a) that there is a unique solution for $\alpha_i$ in $F^3(n) \cap \langle x_n \rangle$. The uniqueness implies that $\alpha_i \in F_2^3(n)$ since, as can be easily checked, $\tau \in F_3^4(n)$. □

We can now use part (*b*) of lemma 2.3 to define a homomorphism $\overline{\tau} : F^3(n-1) \to A(F^4(n))$ as follows. $\overline{\tau}(w)$ is the unique element $\alpha \in A(F^4(n))$ satisfying $\alpha(x_i) \equiv x_i \mod \langle x_n \rangle$ for $i < n$ and $\alpha(x_n) = \tilde{w} x_n \tilde{w}^{-1}$ for some $\tilde{w} \equiv w \mod \langle x_n \rangle$.

Combining this with the restriction $A(F^4(n)) \to A(F^4(n-1))$ we have the following result

**Theorem 1.** *The following sequence is split exact:*

$$1 \longrightarrow F^3(n-1) \xrightarrow{\overline{\tau}} A(F^4(n)) \longrightarrow A(F^4(n-1)) \longrightarrow 1$$

*Remark* 2.5. This proposition is *false* if we replace $F^4(n), F^4(n-1)$ and $F^3(n-1)$ by $F^q(n), F^q(n-1)$ and $F^{q-1}(n-1)$, respectively, if $q > 4$. For example consider the automorphism $\alpha$ of $F^5(2)$ defined by:

$$\alpha(x_1) = [x_2, [x_1, x_2]]x_1[x_2, [x_1, x_2]]^{-1} \tag{11}$$

$$\alpha(x_2) = [x_1, [x_2, x_1]]x_2[x_1, [x_2, x_1]]^{-1} \tag{12}$$

*Proof of theorem* 1. This is an immediate consequence of part (*b*) of lemma 2.3. □

Our next goal is to define a homomorphism $SL^{\mathbb{SE}}(n) \to A(F^4(n))$. First recall the definition of the $\overline{\mu}$-invariants of string links as formulated, for example, in [Le2]. Let $\sigma$ be an n-string link; we define certain canonical *meridian* elements of $\pi(\sigma) = \pi_1(I \times D^2 - \sigma)$. Let $m_i$ be a small circle in $0 \times D^2$ with $0 \times p_i$ as center. Let $s_i$ be a straight line from a base point $p$ in $0 \times (D^2 - \{p_i\})$ to $m_i$. Now let $\mu_i$ be the element of $\pi(\sigma)$ represented by $s_i m_i s_i^{-1}$. Let $\mu : F(n) \to \pi(\sigma)$ be defined by $\mu(x_i) = \mu_i$. It follows from a theorem of Stallings [Sta] that $\mu$ induces an isomorphism $F^q(n) \to \pi^q(\sigma)$ for each $q \geq 3$. We can also define canonical *longitude* elements of $\pi(\sigma)$. Let $l_i$ be a curve on the boundary of a tubular neighborhood of $\sigma_i$ which runs parallel from a point on $m_i$ to the corresponding point on $m_i'$, where $m_i'$ is the projection of $m_i$ into $1 \times D^2$. Let $s_i'$ be the projection of $s_i$ into $1 \times D^2$ and let $u = I \times p$ oriented from 1 to 0. Define $\tilde{\lambda}_i \in \pi$ to be represented by $s_i l_i s_i^{-1} u$. See figure 2. Requiring that $\tilde{\lambda}_i$ have linking number 0 with $\sigma_i$ determines $l_i$ and, hence, $\tilde{\lambda}_i$. We now define $\lambda_i(\sigma) \in F^q(\sigma)$ to be the element (which is $i$- reduced) corresponding to $\tilde{\lambda}_i$ under the isomorphism induced by $\mu$. These elements are invariants of the concordance class of $\sigma$. For a reference on concordance class of string links see [Li] and references therein.

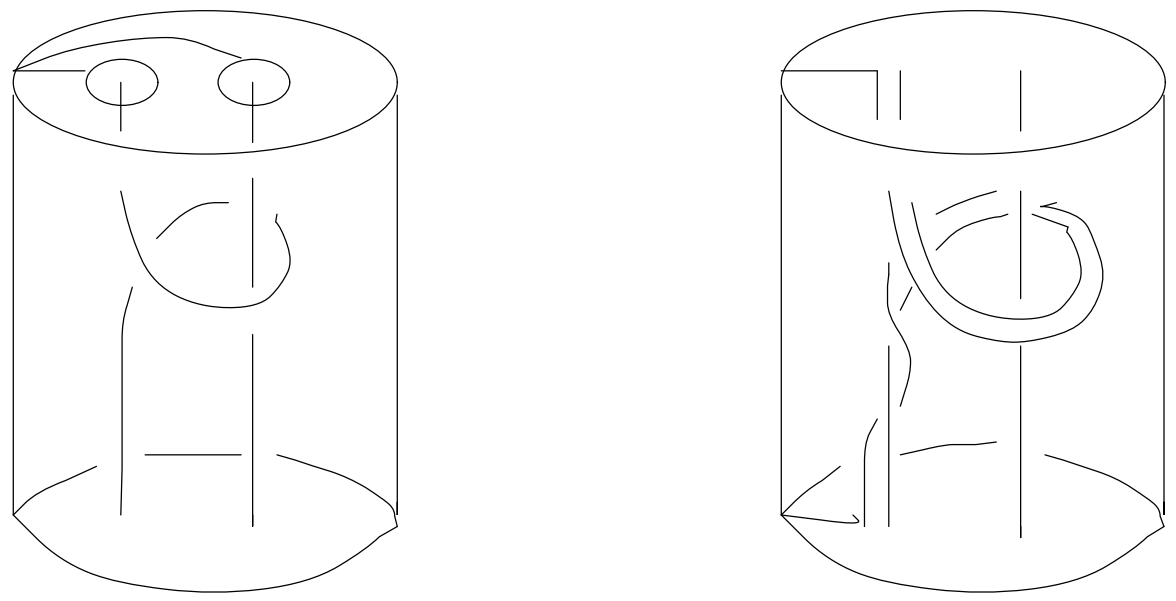

**Figure 2.** On the left shown elements $m_i, s_i$ (for $i = 1, 2$) and on the right shown the element $l_1$ of a string link of 3 components

**Proposition 2.6.** *The surgical equivalence class of $\sigma$ determines the set $\{\lambda_i(\sigma)\} \subset F^3(n)$.*

*Proof.* This is essentially proved in [Le1] but for completeness we give a more direct argument here. Let $\sigma_1$ be an $n$-strand string link and $\gamma$ a simple closed curve in the complement of $\sigma_1$. We assume that the linking numbers of any two strands of $\sigma_1$ and any strand of $\sigma_1$ with $\gamma$ is zero. Let $X = (I \times D^2) - \sigma - \gamma$ and let

$$\pi = \pi_1(X), \quad \pi(\sigma_1) = \pi_1((I \times D^2) - \sigma_1), \quad \pi(\sigma_2) = \pi_1(\Delta - \sigma_1) = \pi_1((I \times D^2) - \sigma_2)$$

If $m$ is a meridian element and $l$ a longitude element for $\gamma$ in $\pi$, then we have:

$$\pi(\sigma_1) \cong \pi/\langle m \rangle \quad \text{and} \quad \pi(\sigma_2) \cong \pi/\langle ml^{\pm 1} \rangle$$

We define homomorphisms:

$$F(n+1) \xrightarrow{\hat{\mu}_1} \pi, \quad F(n+1) \xrightarrow{\hat{\mu}_2} \pi, \quad F(n) \xrightarrow{\mu_1} \pi(\sigma_1), \quad F(n) \xrightarrow{\mu_2} \pi(\sigma_2)$$

all of which send $x_i$ to a meridian of $\sigma_i$ , for $i \leq n$, and, in addition, $\hat{\mu}_1(x_{n+1}) = m$, $\hat{\mu}_2(x_{n+1}) = ml^{\pm 1}$. We have the following commutative diagram:

$$\begin{array}{ccccc} \pi & \xleftarrow{\hat{\mu}_2} & F(n+1) & \xrightarrow{\hat{\mu}_1} & \pi \\ \downarrow & & \downarrow & & \downarrow \\ \pi(\sigma_2) & \xleftarrow{\mu_2} & F(n) & \xrightarrow{\mu_1} & \pi(\sigma_1) \end{array}$$

Here $F(n+1) \to F(n)$ is the obvious reduction defined by sending $x_{n+1}$ to 1. Using results of Milnor [Mi1], [Mi2] and Stallings [Sta] we see that $\mu_1$ and $\mu_2$ induce isomorphisms $\mu_1^3 : F^3(n) \xrightarrow{\cong} \pi^3(\sigma_1)$ and $\mu_2^3 : F^3(n) \xrightarrow{\cong} \pi^3(\sigma_2)$, respectively, while $\hat{\mu}_1$ and $\hat{\mu}_2$ are epimorphisms. Furthermore, there are elements $\lambda_1$ and $\lambda_2$ in $F^3(n+1)$ such that the kernels of $\hat{\mu}_1$ and $\hat{\mu}_2$ (in $F^3(n+1)$) are normally generated, respectively, by $[x_{n+1}, \lambda_1]$ and $[x_{n+1}, \lambda_2]$. In particular, since $\lambda_1, \lambda_2 \in F_2(n+1)$, $\mu_1$ and $\mu_2$ induce isomorphisms $F^3(n) \xrightarrow{\cong} \pi^3$. If $\lambda_{i1}, \lambda_{i2} \in F^3(n)$ are chosen to map to the longitude of $\sigma_i$ in $\pi^3(\sigma_1)$ and $\pi^3(\sigma_2)$, respectively, then they reduce in $F^3(n)$ to $\lambda_i(\sigma_1)$ and $\lambda_i(\sigma_2)$.

Since $\lambda_{i1}, \lambda_{i2} \in F_2^3(n)$, we only need prove that $\hat{\mu}_1$ and $\hat{\mu}_2$ agree on $F_2^3(n+1)$. But this will follow if $\hat{\mu}_2 \equiv \hat{\mu}_1 \mod \pi_2^3$, which is clear since $l \in \pi_2$. □

Now, we can introduce the following definition:

**Definition 2.7.** Let $\Phi_n : SL^{\mathbb{SE}}(n) \to A(F^4(n))$ be the homomorphism defined by

$$\Phi_n(\sigma)(x_i) = \lambda_i(\sigma) x_i \lambda_i(\sigma)^{-1} \tag{13}$$

To see that this is a homomorphism we give another interpretation of $\Phi_n$. Let $X = I \times \partial D^2 \cup 1 \times (D^2 - \{p_i\}) \subseteq I \times D^2$. Then we have the following homomorphisms induced by inclusion maps:

$$\pi_1(0 \times (D^2 - \{p_i\})) \xrightarrow{\mu} \pi_1((I \times D^2) - \sigma) \xleftarrow{\mu'} \pi_1(X)$$

We identify the first and last groups with $F(n)$ equating $x_i$ with the homotopy class of $s_i m_i s_i^{-1}$ and $u^{-1} s_i' m_i' (s_i')^{-1} u$, respectively. Now Stallings theorem tells us that the maps in equation (2.7) become isomorphisms when we pass to any lower central series quotient. Thus the composition of the maps in (2.7) define an automorphism of $F^3(n)$; we leave it as an exercise for the reader to show that this is $\Phi_n(\sigma)$. From this formulation, it is clear that $\Phi_n$ is a homomorphism.

### 2.3. Determining the groups $P^{\mathbb{SE}}(n)$ and $SL^{\mathbb{SE}}(n)$.

In this section we determine the groups $P^{\mathbb{SE}}(n)$ and $SL^{\mathbb{SE}}(n)$.

**Theorem 2.** (a) *$P^{\mathbb{SE}}(n) \to SL^{\mathbb{SE}}(n)$ is an isomorphism.*

(b) *The diagram:*

$$\begin{array}{ccccccccc}
1 & \longrightarrow & F^3(n-1) & \xrightarrow{\tau} & P^{\mathbb{SE}}(n) & \longrightarrow & P^{\mathbb{SE}}(n-1) & \longrightarrow & 1 \\
 & & \| & & \downarrow & & \downarrow & & \\
1 & \longrightarrow & F^3(n-1) & \xrightarrow{\tau'} & SL^{\mathbb{SE}}(n) & \longrightarrow & SL^{\mathbb{SE}}(n-1) & \longrightarrow & 1
\end{array}$$

*is commutative and the rows are split short exact.*

The exactness and commutativity in $(b)$ is clear except for the injectivity of $\tau$ and $\tau'$. This will be proved in theorem 3 below. Clearly $(a)$ will follow from $(b)$.

**Theorem 3.** $\Phi_n$ *is an isomorphism.*

*Proof of theorems* 2 *and* 3. We can combine the statements of theorems 1 2 and 3 into a single commutative diagram:

$$\begin{array}{ccccccccc}
1 & \longrightarrow & F^3(n-1) & \xrightarrow{\tau} & P^{\mathbb{SE}}(n) & \longrightarrow & P^{\mathbb{SE}}(n-1) & \longrightarrow & 1 \\
 & & \| & & \downarrow & & \downarrow & & \\
1 & \longrightarrow & F^3(n-1) & \xrightarrow{\tau'} & SL^{\mathbb{SE}}(n) & \longrightarrow & SL^{\mathbb{SE}}(n-1) & \longrightarrow & 1 \\
 & & \| & & \downarrow & & \downarrow & & \\
1 & \longrightarrow & F^3(n-1) & \xrightarrow{\overline{\tau}} & A_0(F^4(n)) & \longrightarrow & A_0(F^4(n)) & \longrightarrow & 1
\end{array}$$

The commutativity of this diagram is clear. To establish that the first two rows are exact we need only confirm the injectivity of $\tau$ and $\tau'$. But this now follows from the injectivity of $\overline{\tau}$. This concludes the proofs of theorems 2 and 3. □

**Corollary 2.8.** $P^{\mathbb{SE}}(n)$ *is a nilpotent group of class two, i.e., with the notation of section* 2.1 *we have that* $P_3^{\mathbb{SE}}(n) = 1$. *Furthermore,* $P^{\mathbb{SE}}(n)/P_2^{\mathbb{SE}}(n)$ *is a free abelian group in* $\binom{n}{2}$ *generators and* $P_2^{\mathbb{SE}}(n)/P_3^{\mathbb{SE}}(n)$ *is a free abelian group in* $\binom{n}{3}$ *generators. In fact,* $P_2^{\mathbb{SE}}(n)/P_3^{\mathbb{SE}}(n)$ *is a free abelian group on generators* $\alpha_I$ *for* $I \in I_n := \{(i,j,k)|1 \le i,j,k \le n\}$ *(where all* $i,j,k$ *are distinct) with identifications* $\alpha_I = sgn(\sigma)\alpha_{\sigma(I)}$, *where* $\sigma$ *is any permutation of* $I$. *The element* $\alpha_I$ *can be represented, if* $i < j < k$, *by the braid* $B_I$ *in figure* 3.

*Remark* 2.9. The relation of surgical equivalence on string links and pure braids is generated by the local moves of figures 4 and 5.

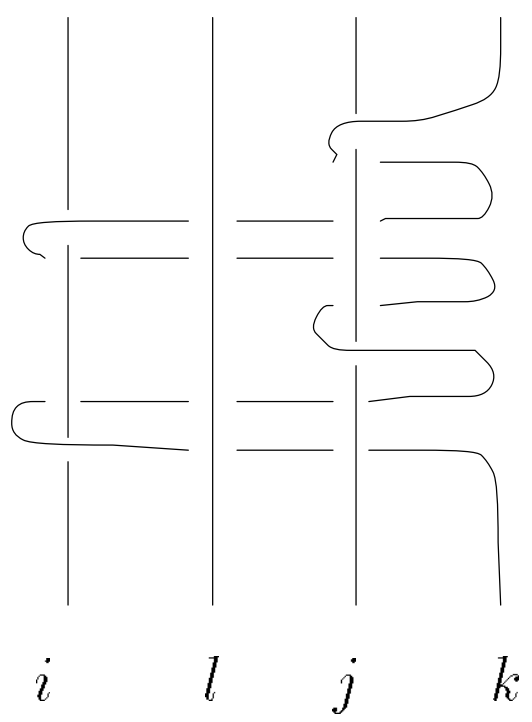

**Figure 3.** A 4 strand braid $B_{i,j,k} = B_1 \cup B_2$ representing $\alpha_{i,j,k}$. Here $B_1$ consists of the 3 strands $i, j, k$ and $B_2$ consists of all the rest strands (here only one). Note that the strands of $B_2$ are on top of the $i, j, k$ strands, and that the closure of the 3 strand braid $B_2$ is a Borromean link of 3 components and that the link represented by the closure of the $i, j$ strands is trivial. Furthermore, the longitude $l_k$ of the $k^{th}$ strand represents the element $[\mu_i, \mu_j]$ in the fundamental group of the complement of the link consisting of the $i, j$ strands, where $\{\mu_j\}$ are the canonical free generators of the fundamental group of the complement of the $i, j$ strands.

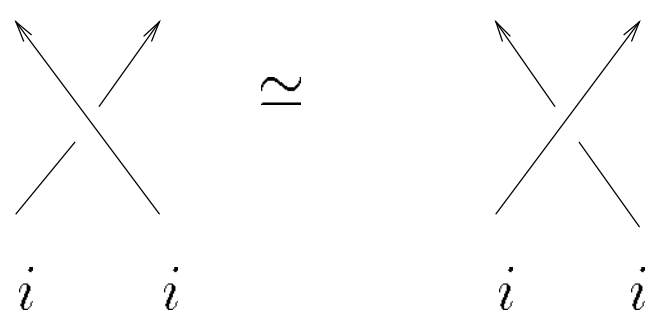

**Figure 4.** A local move that generates the equivalence relation of string link homotopy. Here $i$ denotes the same component of a string link.

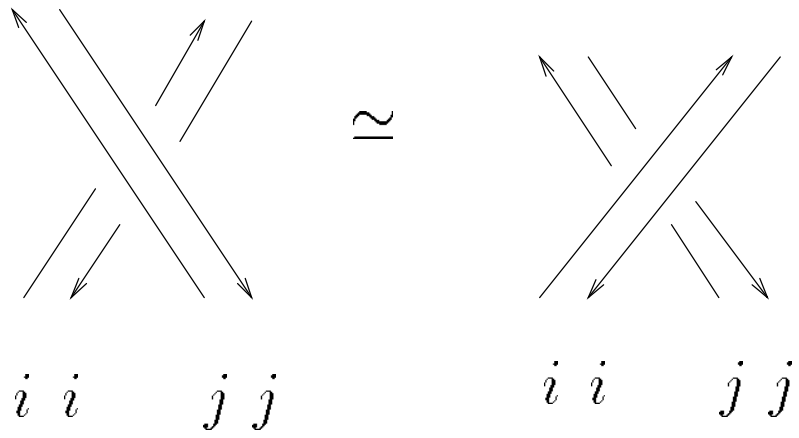

**Figure 5.** A local move that is implied by the relation of surgical equivalence of string links. Here $i, j$ denote different components of a string link.

## 3. Relations with finite type invariants of $\mathbb{Z}HS$

**3.1. Preliminaries on finite type 3-manifold invariants.** As a motivation to the notion of finite type 3-manifold invariants, let us recall the definition of type $m$ Vassiliev invariants of (oriented) knots in (oriented) $S^3$ (after [B-N], [BL], [Va]): $V$ is type $m$ invariant of knots if for every knot $K$ in $S^3$ and every choice $S_1, \dots, S_{m+1}$ of embedded 2-spheres that intersect the knot as in figure 6 we have that

$$\sum_{\epsilon_i \in \{0,1\}} \prod_{i=1}^{m+1} (-1)^{\epsilon_i} V(K_{\epsilon_1, \dots, \epsilon_{m+1}}) = 0 \tag{14}$$

where $K_{\epsilon_1, \dots, \epsilon_{m+1}}$ is the knot obtained by cutting $S^3$ along $S_1, \dots, S_{m+1}$, twisting every $S_i$ $\epsilon_i$ times as in figure 6 and gluing back.

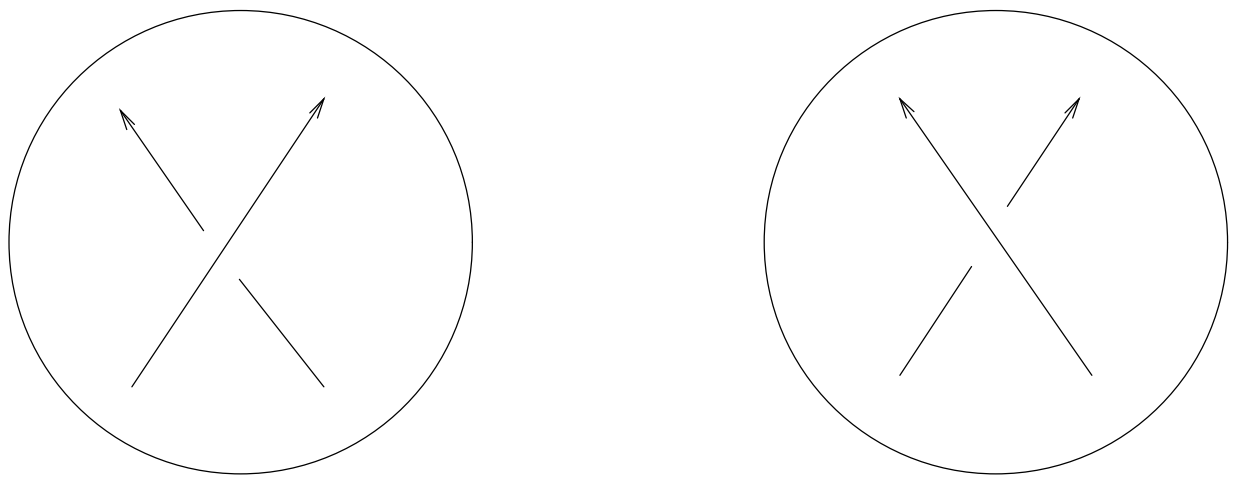

**Figure 6.** An embedded 2-sphere in $S^3$ intersecting a knot in as in the left part of the figure. Also shown on the left is the result of an $\epsilon = 0$ diffeomorphism (i.e. no twist, in other words the identity diffeomorphism) of the 2-sphere and and on the right is the result of an $\epsilon = 1$ diffeomorphism of the 2-sphere (i.e. a $\pi$ twist along the vertical axis of the page you are reading).

Let $\mathcal{F}_m\mathcal{V}$ be the space of type $m$ knot invariants, and let $\mathcal{F}_\star\mathcal{V} := \cup_{m\geq 0}\mathcal{F}_m\mathcal{V}$ be the space of finite type knot invariants. It is easy to show that $\mathcal{F}_\star\mathcal{V}$ is a filtered commutative algebra (with pointwise multiplication). Let $\mathcal{G}_\star\mathcal{V}$ denote the associated graded algebra (and more generally, let $\mathcal{G}_\star(obj)$ denote the associated graded object of the filtered object $\mathcal{F}_\star(obj)$) and let $\mathcal{W}_m$ denote the space of weight systems of degree $m$ (i.e. linear functionals on the space of chord diagrams with $m$ chords modulo $4T$ relations and framing independence relation, see references above). Then it is easy to show that there is a map $\mathcal{F}_m\mathcal{V} \to \mathcal{W}_m$ with kernel $\mathcal{F}_{m-1}\mathcal{V}$. Since $\mathcal{W}_m$ is *a priori* finite dimensional, so is $\mathcal{F}_m\mathcal{V}$.

Before we talk about 3-manifold invariants, let us establish some useful notation: Let $\mathcal{M}$ denote the vector space (over $\mathbb{Q}$) on the set of oriented integral homology 3-spheres ($\mathbb{Z}HS$ for short). A link $L \subseteq M$ in a $\mathbb{Z}HS$ is called algebraically split if the linking numbers between any two components vanish. A framing $f = (f_1, \dots, f_n)$ for an $n$ component link is a sequence of integers associated to each component.

*Remark* 3.1. A framing for a link $L$ corresponds to a choice of *longitudes* for each component. If $N_i$ is a tubular neighborhood of the component $L_i$, then the associated longitude $\gamma_i$ on $\partial N_i$ is required to be homologous to $L_i$ in $N_i$ and to have linking number $f_i$ with $L_i$. To make sense of this we must impose an orientation on $L_i$, but then it is easy to see that $\gamma_i$ is independent of this choice. Note that linking numbers make sense in any ZHS.

A framed link $(L, f)$ in a $\mathbb{Z}HS$ $M$ is called unimodular if $f_i = \pm 1$ for all $i$. A framed link $(L, f)$ is called $AS$ admissible if it is algebraically split and unimodular. For every framed link $(L, f)$ in $M$ we denote by $M^{L,f}$ the result of doing *Dehn surgery* on $(L, f)$ in $M$ [Ro], i.e. remove a tubular neighborhood $N_i$ of each $L_i$ and sew it back in so that $\gamma_i$ now bounds a disk in $N_i$. For a framed link $(L, f)$ in $M$ we denote

$$(M, L, f) := \sum_{L' \subseteq L} (-1)^{|L'|} M^{L', f|_{L'}} \in \mathcal{M} \tag{15}$$

Recall that $\mathcal{M}$ is the $\mathbb{Q}$ vector space on the set of $\mathbb{Z}HS$, and that if $(L, f)$ is a framed link in a $\mathbb{Z}HS$ $M$, then $(L, f)$ is $AS$ admissible if and only if $M^{L', f|_{L'}}$ is a $\mathbb{Z}HS$ for every sublink $L'$ of $L$. Let us define a decreasing filtration $\mathcal{F}^{Oh}_\star$ on the vector space $\mathcal{M}$ as follows: $\mathcal{F}^{Oh}_m \mathcal{M}$ is the subspace spanned by $(M, L, f)$ for all $AS$ admissible unimodular links of $m$ components.

We can equally well consider $\mathcal{M}$ to be the vector space generated by all triples $(M, L, f)$ subject to the *fundamental relation*:

$$(M, L, f) = (M, L', f|L') - (M^{(l, f|_l)}, L', f|L') \tag{16}$$

where $L'$ is any sublink of $L$ obtained by removing a component $l$. If $l$ bounds a disk $D$ in $M$, then $M^{(l, f|_l)} \cong M$ and we may construct the link in $M$ corresponding to $L$ in $M^{(l, f|_l)}$ from $L$ by just giving the bundle of strands of $L$ which pass through $D$ a full clockwise twist, as in figure 7.

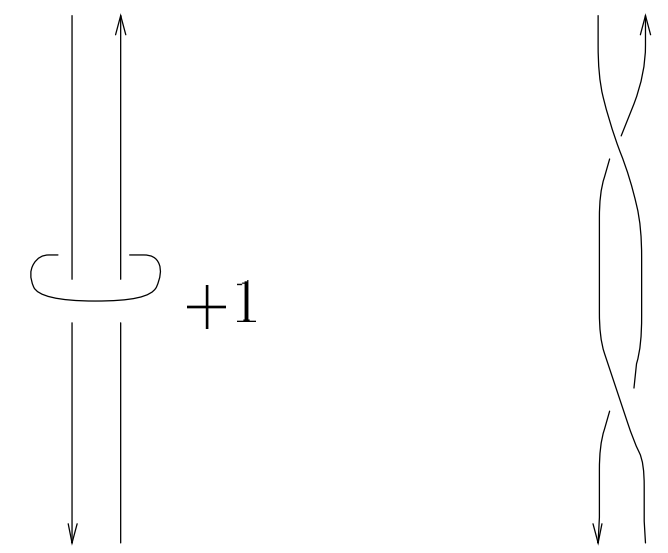

**Figure 7.** A framed link of two components on the left (only two components are shown) represents (after Dehn surgery) the same 3-manifold as the link of one component less on the right. The link on the right is obtained by a $-1$ twist.

We can now recall the following definition from [Oh]:

**Definition 3.2.** $\lambda$ is a type $m$ invariant of $\mathbb{Z}HS$ (with values in $\mathbb{Q}$) if $\lambda(\mathcal{F}^{Oh}_{m+1}\mathcal{M}) = 0$ i.e., if for every $AS$ admissible link $L$ of $m+1$ components in a $\mathbb{Z}HS$ $M$ we have that

$$\sum_{L' \subseteq L} (-1)^{|L'|} \lambda(M^{L', f|_{L'}}) = 0 \tag{17}$$

Let $\mathcal{F}_m\mathcal{O}$ denote the vector space of type $m$ invariants of $\mathbb{Z}HS$, and let $\mathcal{F}_\star\mathcal{O}$ denote the union $\cup_{m\geq 0}\mathcal{F}_m\mathcal{O}$. It is easy to see that $\mathcal{F}_\star\mathcal{O}$ is a filtered commutative algebra with pointwise multiplication.

**3.2. Surgical equivalence and 3-manifold invariants.** In this section we reprove Ohtsuki's fundamental result (theorem 4) which states that the space of type $m$ invariants of $\mathbb{Z}HS$ is finite dimensional for every $m$.

We begin by observing that every $\sigma \in P^{\mathbb{I}}(n)$ can be closed to a link $\hat{\sigma}$ of $n$ components in $S^3$. Furthermore, with the notation of corollary 2.8 we have that $\sigma \in P^{\mathbb{I}}_2(n)$ if and only if $\hat{\sigma}$ is an algebraically split link. Let us now consider the map $P^{\mathbb{I}}_2(n) \to \mathcal{M}$ defined as $\sigma \to (S^3, \hat{\sigma}, \{+1, \dots, +1\})$. We claim that this map descends to a well defined map ( *not* a group homomorphism)

$$P^{\mathbb{SE}}_2(n) \longrightarrow \mathcal{G}^{Oh}_n\mathcal{M} \tag{18}$$

Indeed, it follows from the definition of $\mathbb{SE}$ equivalence of pure braids. This extends to a linear map from the rational group ring $\mathbb{Q}[P^{\mathbb{SE}}_2(n)] \to \mathcal{G}^{Oh}_n\mathcal{M}$. By corollary 2.8 we can regard $P^{\mathbb{SE}}_2(n)$ as a free abelian group on generators $\alpha_I$ for $I \in I_n := \{(i,j,k) | 1 \leq i,j,k \leq n\}$ (where all $i,j,k$ are distinct) with identifications $\alpha_I = \text{sgn}(\sigma)\alpha_{\sigma(I)}$, where $\sigma$ is any permutation of $I$. Therefore, $\mathbb{Q}[P^{\mathbb{SE}}_2(n)]$ is the Laurent polynomial ring in the commuting variables $\alpha_I$.

*Remark* 3.3. Using the fundamental relation (16) we see that if an $n$ $AS$ admissible link $L'$ is obtained by twisting along an admissible unknot $\gamma$ an $AS$ admissible $n$ component link $L$, then we have that

$$(S^3, L \cup \gamma, \{+1, \dots, +1\}) = (S^3, L\{+1, \dots, +1\}) - (S^3, L', \{+1, \dots, +1\}) \tag{19}$$

Therefore, surgical equivalence relates $n$ component $AS$ admissible links in $S^3$ modulo $n+1$ component $AS$ admissible ones. Since $P^{\mathbb{SE}}(n)$ is not a finite group (instead, it is a finitely generated free abelian group) at this point it is not clear why $\mathcal{G}^{Oh}_n\mathcal{M}$ is a finite dimensional $\mathbb{Q}$ vector space.

Our next task is to introduce a finite dimensional quotient of $\mathbb{Q}[P^{\mathbb{SE}}_2(n)]$ that maps onto $\mathcal{G}^{Oh}_n\mathcal{M}$. This will involve looking at $n$ and $n-1$ component $AS$ admissible links in $S^3$. First, we define a monomial $a = \prod_{i \in I_n} \alpha_I^{a(I)}$ to be $i$-trivial (for $1 \leq i \leq n$) if $a(I) \neq 0$ implies $i \in I$ and to be $i$-disjoint if $a(I) \neq 0$ implies $i \notin I$. Now we define our set of relations to be all those obtained by taking linear combinations with rational coefficients of the following elementary relations:

$$\begin{aligned}
(20)\quad && a = & \quad 0 && \text{if } a \text{ is } i\text{-disjoint for some } i \\
(21)\quad && d(a-1)(b-1)(c-1) = & \quad -d && \text{if } a, b, c \text{ are } i\text{-trivial and} \\
&& & && d \text{ is } i\text{- disjoint for some } i \\
(22)\quad && d(a-a^{-1}) = & \quad 0 && \text{if } a \text{ is } i\text{-trivial and} \\
&& & && d \text{ is } i\text{-disjoint for some } i
\end{aligned}$$

We will denote by $\mathcal{R}_n$ the subspace of $\mathbb{Q}[P_2^{\mathbb{SE}}(n)]$ generated by equations (20), (21) and (23).

Some useful particular consequences of these relations are:

$$\begin{aligned}
(23)\quad && 1 = & \quad 0 && \\
(24)\quad && da^2 = & \quad d(4a-2) && \text{if } a \text{ is } i\text{-trivial} \\
&& & && \text{and } d \text{ is } i\text{-disjoint for some } i \\
(25)\quad && d(ab^{-1}+ab) = & \quad d(2a+2b-1) && \text{if } a, b \text{ are } i\text{-trivial and} \\
&& & && d \text{ is } i\text{-disjoint for some } i
\end{aligned}$$

Indeed, setting Setting $a = 1$ in (20) gives (23). Setting $b = a, c = a^{-1}$ in (21) and using (20) gives (24). Setting $c = b^{-1}$ in (21) and using (20) gives (25). □

**Theorem 4 (Oh).** *The map $\mathbb{Q}[P_2^{\mathbb{SE}}(n)] \to \mathcal{G}_n^{Oh}\mathcal{M}$ is onto. Furthermore, it factors through a (necessarily onto) map $\mathbb{Q}[P_2^{\mathbb{SE}}(n)]/\mathcal{R}_n \to \mathcal{G}_n^{Oh}\mathcal{M}$.*

*Proof of theorem* 4. We first show ontoness. Recall that $\mathcal{G}_n^{Oh}\mathcal{M}$ is generated by triples $(M, L, f)$ where $L$ is an $AS$ admissible link of $n$ components in a $\mathbb{Z}HS$ $M$. Every 3-manifold can be obtained from $S^3$ by surgery on a framed link and, if the manifold is a $\mathbb{Z}HS$, this means that the matrix of linking numbers of the link is unimodular. By a sequence of handle slides we can diagonalize the linking matrix and so we can assume that $M$ is obtained by surgery on an $AS$ admissible link in $S^3$. Using our fundamental relation (16) we can assume that $M = S^3$. Ontoness now follows once we can show that we can assume that the framing $f_i$ is equal to $+1$, for all $1 \leq i \leq n$. This follows by an argument due to Ohtsuki [Oh]. Let $L$ be any $AS$ admissible $n$ component link with two framings $f$ and $f'$ so that, if $L'$ is the sublink obtained by deleting the $n$-th component, then $f|L' = f'|L'$ and $f_n = -f'_n$. Let $(\tilde{L}, \tilde{f})$ be obtained from $(L, f)$ by replacing $L_n$ by two parallel non-linking copies of $L_n$, and with $\tilde{f}|L' = f|L'$, $\tilde{f}_n = f_n$ and $\tilde{f}_{n+1} = -f_n$. Then we have $(M, \tilde{L}, \tilde{f}) = (M, L, f) - (M^-, L, f)$ where $M^\epsilon = M^{(L_n, \epsilon f_n)}$. So it suffices to show that $(M^-, L, f) = -(M, L, f')$. But we have $(M, L, f') = (M, L', f'|L') - (M^-, L', f'|L')$ and $(M^-, L, f) = (M^-, L', f'|L') - (M^{-+}, L', f'|L')$. Since $M^{-+} = M$ and $f|L' = f'|L'$, we are done. The same proof is represented in the following figure 8:

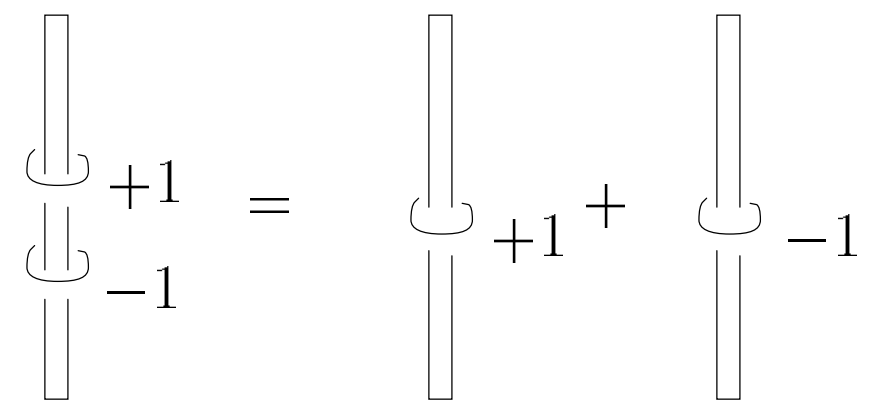

**Figure 8.** A proof that $+1$ framings in each component suffice to generate $\mathcal{G}_n^{Oh}\mathcal{M}$

To prove that the map $\mathbb{Q}[P_2^{\mathbb{SE}}(n)] \to \mathcal{G}_n^{Oh}\mathcal{M}$ factors through a map $\mathbb{Q}[P_2^{\mathbb{SE}}(n)]/\mathcal{R}_n \to \mathcal{G}_n^{Oh}\mathcal{M}$ we recast the relations in $\mathcal{R}_n$ in terms of links.

To prove equation (20), we map $a$ to $(S^3, L, f)$ where the $i$- th component $l$ of $L$ bounds a disk disjoint from the remaining components $L'$. Since surgery on $l$ does not change $L'$ this relation follows from the fundamental relation (16).

To prove equation (21), after closing up the braid $abcd$ to a link $L$, the $i$-th component $l$ represents the element $\alpha\beta\gamma$ in the fundamental group of the complement of the remaining components $L'$ and the surgical equivalence class of $L$ only depends on this element. Furthermore we can assume that $l$ is unknotted since this can be achieved by crossing changes which means surgeries along small circles enclosing the crossings and the class of $L$ in $\mathcal{G}_n^{Oh}\mathcal{M}$ is unchanged by surgeries. In this way we can reformulate equation (21) and, by a similar argument, equation (23), as follows.

Let $L$ be a $(m-1)$-component link in $S^3$ and $\alpha$ be an element of $\pi = \pi_1(S^3 - L)$. Then we can consider the $m$-component link $L(\alpha)$ defined by adding to $L$ a new component which represents $\alpha$. Note that $L(\alpha)$ is well-defined up to homotopy and, therefore, up to surgical equivalence. $L(\alpha)$ is algebraically split if and only if $L$ is algebraically split and $\alpha \in \pi_2$, the commutator subgroup of $\pi$.

**Lemma 3.4.** *Suppose that $L$ is an algebraically split $(m-1)$- component link and $\alpha_1, \dots, \alpha_n \in \pi_2$ and $n \geq 3$. Then, we have the following identities in $\mathcal{G}_m^{Oh}\mathcal{M}$:*

$$L(\alpha_1 \dots \alpha_n) = \sum_{1 \leq i < j \leq n} L(\alpha_i\alpha_j) - (n-2) \sum_{1 \leq i \leq n} L(\alpha_i) \tag{26}$$

$$L(\alpha^{-1}) = L(\alpha) \tag{27}$$

*Proof.* Equation (27) follows from the fact that surgery along a simple closed curve is independent of the orientation of the curve. Equation (26) will follow from a more general relation in $\mathcal{M}$ (theorem 5) stated below. □

Before we state theorem 5, we need to fix some notation: Consider a disk $D$ embedded in $S^3$ with $n \geq 3$ bundles of strands, from a link $L$ disjoint from $\partial D$, passing through $D$. Assume the strands from any single component of $L$ in any single bundle have an equal number in each direction. We consider links formed by adjoining to $L$ circles in $D$ enclosing one or more of the bundles as in figure 9. We denote them by $\overline{x_i}, \overline{x_i x_j}, \overline{x_i x_j x_k}$, as indicated in figure 9. Moreover we multipy them by layering the corresponding circles in parallel copies of $D$, as indicated in figure 9. Note that this multiplication is *not* commutative in $\mathcal{M}$, although it is in $\mathcal{G}_n^{Oh}\mathcal{M}$. However $\overline{x_i}$, does commute with everything. In [GO] we will explore this noncommutativity further.

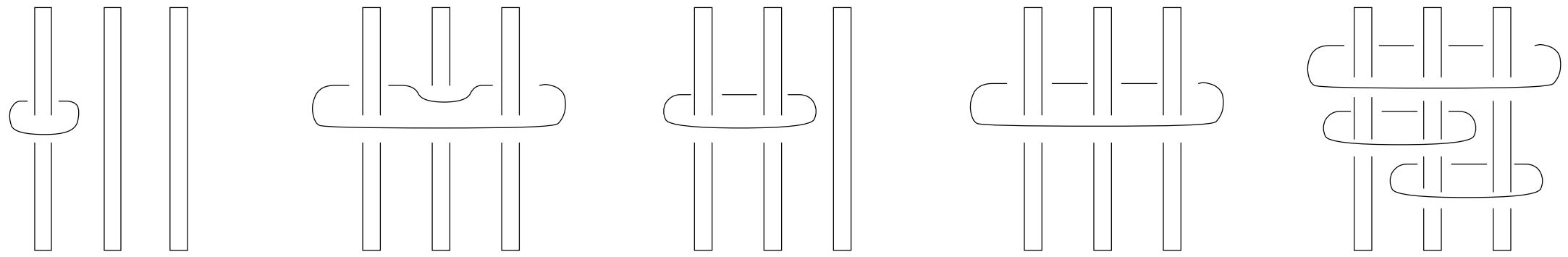

**Figure 9.** Figures representing $\overline{x_1}, \overline{x_1 x_2}, \overline{x_1 x_3}, \overline{x_1 x_2 x_3}$. On the right is shown a composition of figures by stacking from bottom to top.

**Theorem 5.** *With the above notation we have the following:*

$$1 - \overline{x_1 \dots x_n} = \frac{\prod_{1 \le i \le n} \prod_{i < j \le n} (1 - \overline{x_i x_j})}{\prod_{\le i \le n} (1 - \overline{x_i})^{n-2}} \tag{28}$$

*This is really an equation in the completion of* $\mathcal{M}$ *or in* $\mathcal{M}/\mathcal{F}_m^{Oh}\mathcal{M}$ *for any* $m$.

*Proof.* The proof is a refinement of Ohtsuki's idea of resolving a full twist of three bands into a sequence of twists of pairs and individual bands. For simplicity, we assume that $n = 3$. The proof is given in pictures 13 and 14. The drawing conventions are given in pictures 10, 11 and 12. Figure 14 follows from figures 11 and 12. Note that the framings in each component represented by horizontal unknots in figure 13 is $+1$. □

The proof of lemma 3.4, and therefore of theorem 4, is now complete. □

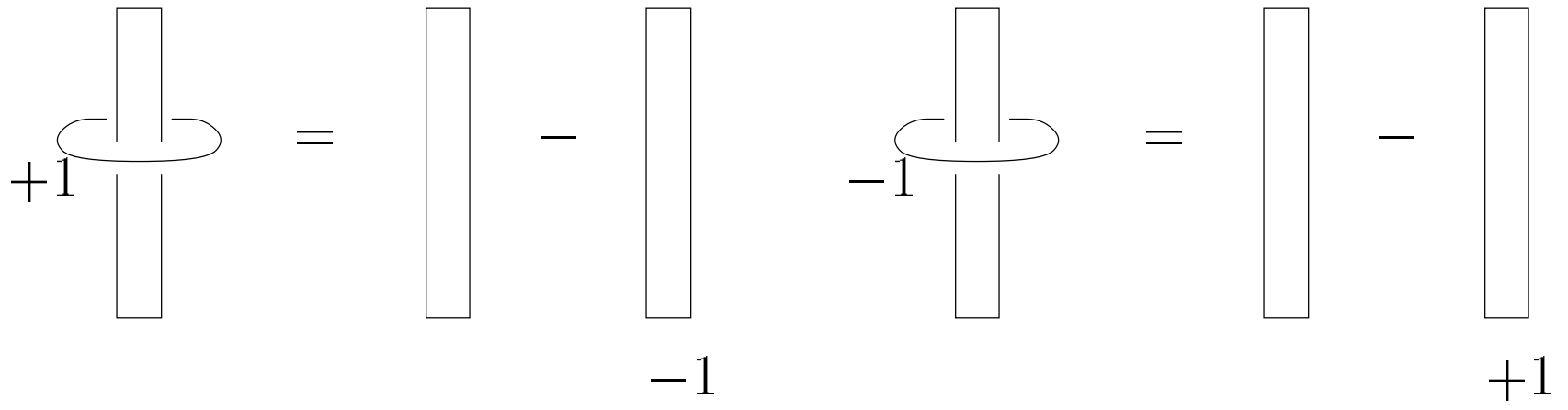

**Figure 10.** Some conventions on drawing $AS$ admissible links in $\mathbb{Z}HS$. Here a band stands for a knot (in the $\mathbb{Z}HS$ represented by $+1$ surgery on an unknot). A (negative or positive) twist on a band is shown on the right.

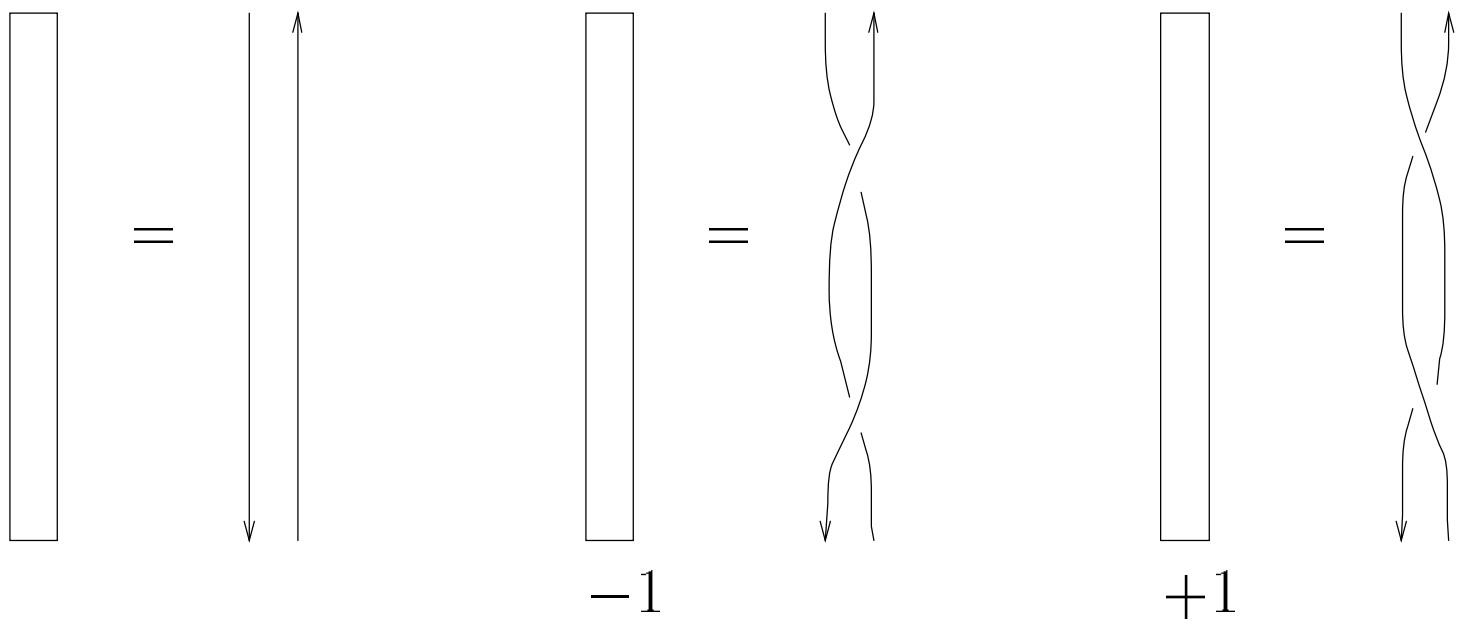

**Figure 11.** Some equalities in $\mathcal{M}$ coming from Kirby calculus. Shown here are a ribbon part of $AS$ admissible links that represent (linear combinations of) $\mathbb{Z}HS$ under the map (15). The numbers in the bottom of each band indicate the number of twists that we put in the band.

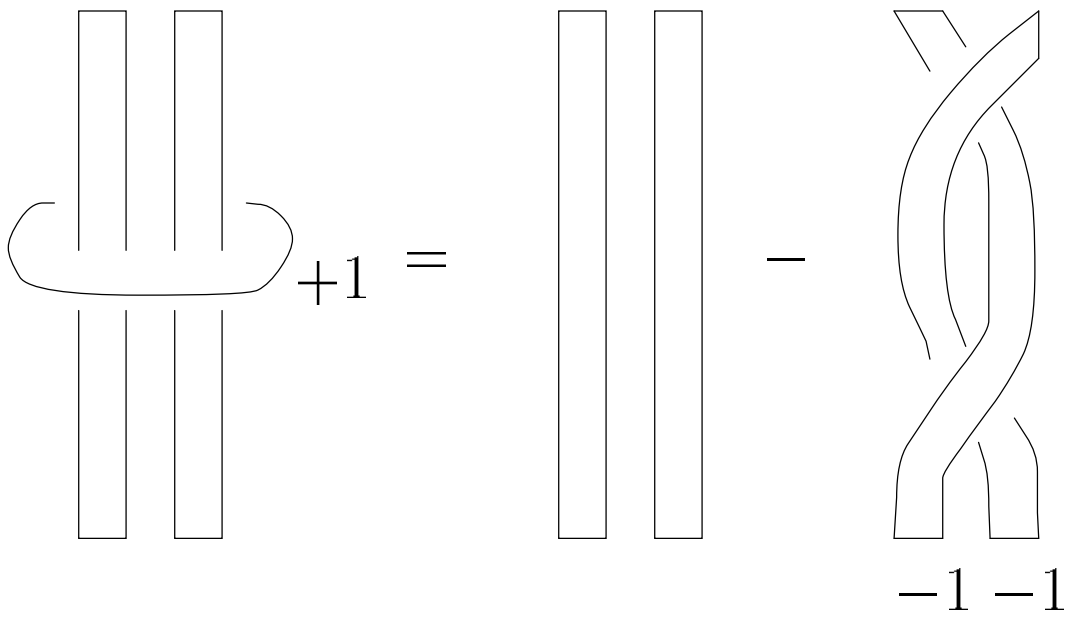

**Figure 12.** Some more equalities in $\mathcal{M}$ coming from Kirby calculus. Shown here are two ribbon parts of $AS$ admissible link that represent (linear combinations of) $\mathbb{Z}HS$ under the map (15). The two ribbons may represent the same component or not.

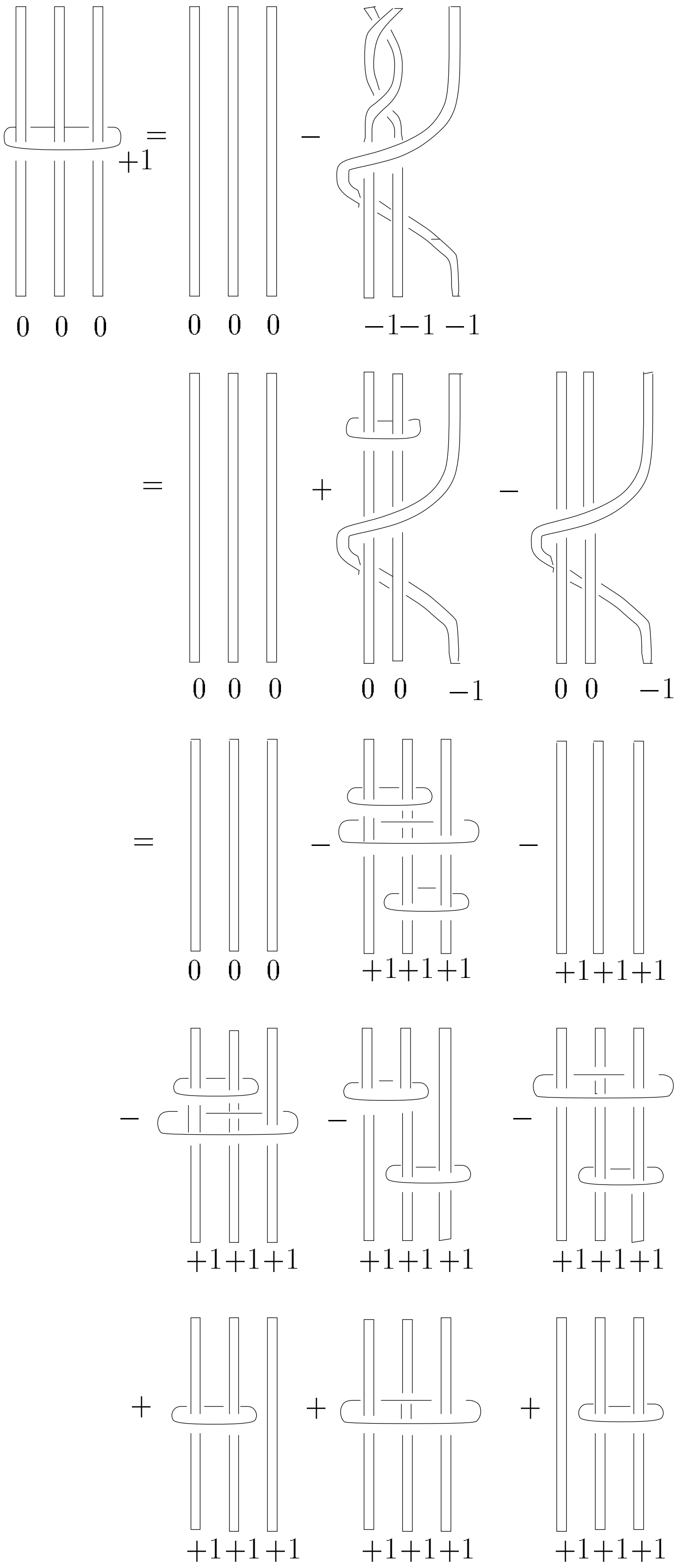

**Figure 13.** The beginning of the proof of theorem 5. The two equalities here follow from repeated applications of figures 11 and 12. The framings in all horizontal components is $+1$.

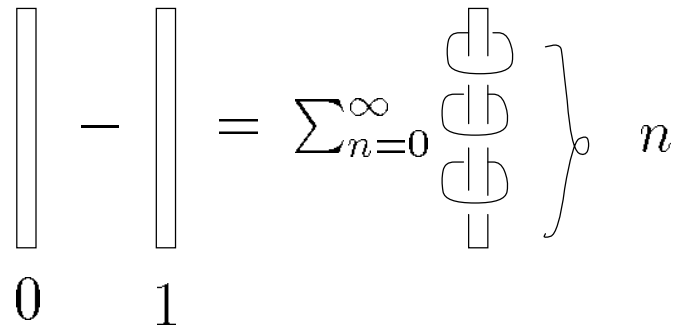

**Figure 14.** The end of the proof of theorem 5. The identity in the present figure follows from figure 8. The framings of the horizontal components represented by unknots is +1.

**3.3. A vanishing theorem for finite type invariants of $\mathbb{Z}HS$.** In this section we prove a vanishing theorem (theorem 6 and corollary 3.5) for the graded vector space of finite type invariants of $\mathbb{Z}HS$. With the notation of the introduction (see sections 1.1, 1.2), recall that $\mathcal{M}$ is the vector space (over $\mathbb{Q}$) on the set of oriented $\mathbb{Z}HS$, and elements in $\mathcal{F}_m\mathcal{O}$ are thought of as polynomial functions of degree $m$ on $\mathcal{M}$.

Consider the following three properties of a graph $G$:

(a) Every edge has two distinct vertices.
(b) Every vertex is either univalent or trivalent.
(c) Every component of $G$ is either a $Y$-graph (i.e. has exactly one trivalent vertex, three univalent vertices and three edges) or has every vertex trivalent.
(d) $G$ contains no $\Theta$ component, i.e. a trivalent graph with 3 edges and 2 vertices.

We will say $G$ is a *UT graph* if it satisfies (a) and (b), and an *Ohtsuki graph* if it also satisfies (c) and (d). The *degree* of $G$ is the number of edges. It is easy to see that any Ohtsuki graph has degree a multiple of 3.

We can now state the following theorem:[2]

**Theorem 6.** *For every $n \geq 0$, the vector space $\mathbb{Q}[P_2^{\mathbb{SE}}(n)]/\mathcal{R}_n$ has a basis in one-one correspondence with the set of Ohtsuki graphs of order $n$.*

*Proof of theorem* 6. For a monomial $a = \prod_{I \in I_n} \alpha_I^{a(I)}$ and for $1 \leq i \leq n$ we define the $i$-order of $a$, $o_i(a)$, to be $\sum_{i \in I} |a(I)|$ and the *total order* of $a$ to be $\sum_I |a(I)|$. We show that $\mathbb{Q}[P_2^{\mathbb{SE}}(n)]/\mathcal{R}_n$ is generated by those $a$ which satisfy: $o_i(a) \leq 2$ for every $i$, and $a(I) = \pm 1$ for every $I$- we call such $a$ *admissible*. If $o_i(a) > 2$ then we may write $a = \alpha_{I_1}^{\epsilon_1} \alpha_{I_2}^{\epsilon_2} bc$, for some $i \in I_1 \cap I_2$, where $b$ is $i$-trivial and $c$ is $i$-disjoint. By equation (21) we see that $a$ can be written as a linear combination of monomials with strictly smaller $i$-order and no larger total order. So we can assume $o_i(a) \leq 2$ for every $i$. Now suppose $a(I) = \pm 2$ for some $I$. Then $a = \alpha_I^{a(I)} b$ where $b$ is $i$-disjoint and $i \in I$. We can apply equation (23) and (24) to express $a$ as a linear combination of monomials with strictly smaller $i$-order and no larger total order.

We can, following Ohtsuki, associate to any admissible monomial $a$ a graph $G(a)$ as follows. For every $i = 1, \ldots, n$ we associate an edge $e(i)$ with two distinct vertices.

[2]A more refined version of theorem 6 will appear in [GO].

For every $I$ with $a(I) \neq 0$ we define a trivalent vertex $v(I)$ incident to each $e(i)$ with $i \in I$. We can also give each trivalent vertex $v(I)$ an *orientation*, i.e. a cyclic ordering to its incident edges, by choosing the ordering given by those $I$ for which $a(I) = +1$. It is easy to see that $G(a)$ is a vertex oriented $UT$ graph with no $\Theta$ components. Furthermore, $G(a)$ determines $a$. We can use this graphical interpretation to describe further reductions to the generating set of monomials. We will say an edge is interior if both its vertices are trivalent and exterior otherwise.

We make the following observations:

(a) If $G(a)$ has an edge with both vertices univalent, which will be true if the order of $a$ is small enough, then $a = 0$; this follows from (20).
(b) Suppose $e$ is any edge in $G(a)$. Define $a'$ so that $G(a) = G(a')$ except that the trivalent vertices of $e$ have orientations reversed. Then (23) implies $a = a'$.
(c) If $e$ is an interior edge and we define $a''$ by changing the orientation of only one of the vertices, then (25) says that $a + a''$ is a linear combination of monomials with smaller total order.

Suppose now that $G(a)$ has two edges, one interior and the other exterior, which share a common vertex. Then (b) and (c) imply that $2a$ is a linear combination of monomials with smaller total order. But if $G(a)$ contains no such edges it is easy to see that it must be an Ohtsuki graph. Thus the set of all such graphs *with oriented vertices* form a generating set. To complete the proof we choose some preferred orientation for each Ohtsuki graph and show that, if $a$ is any monomial corresponding to this graph with *another* orientation then it is equal to $\pm$ the preferred one plus monomials of lower total order. To see this let $G(a)$ be the associated graph and let $G(a')$ be produced by changing the orientation of one of the trivalent vertices $v$. Then $a' = \pm a +$ monomials of lower total order. In fact, if $v$ lies in a Y component then, by (a), $a = a'$ while if $v$ lies in a completely trivalent component then, by (b), $a' = -a +$ monomials of lower total order.

To prove linear independence we reinterpret $\mathbb{Q}[P_2^{\mathbb{SE}}(n)]/\mathcal{R}_n$ as the quotient of the vector space spanned by *unreduced* but commuting monomials in the $\{\alpha_I\}$, i.e. $\alpha_I$ and $\alpha_I^{-1}$ can both occur in a monomial, by the relations defining $\mathcal{R}_n$ and the cancellation relation:

$$aa^{-1}bd = d \quad \text{if } a, b \text{ are } i\text{-trivial and } d \text{ is } i\text{-disjoint for some } i \tag{29}$$

We write an unreduced monomial uniquely in the form $\alpha = \alpha_{I_1}^{e_1} \cdots \alpha_{I_n}^{e_n}$, where each $I_r = ijk$ with $i < j < k$ and $e_r = \pm 1$. Then we define the *multiplicity* $m_i = \sum_{i \in I_r} |e_r|$, for each $i = 1, \ldots, n$. We say $i$ is *multiple* if $m_i > 1$ and *singular* if $m_i = 1$. We now define a canonical reduction of $\alpha$ in five steps.

Step 1. Write:

$$\alpha = \sum \alpha_{I_{i_1}}^{e_{i_1}} \cdots \alpha_{I_{i_k}}^{e_{i_k}} - \sum c_{i_1 \cdots i_k} \alpha_{I_{i_1}}^{e_{i_1}} \cdots \alpha_{I_{i_k}}^{e_{i_k}}$$

where the first summation is over all subsequences of $I_1, \ldots, I_n$ such that the multiplicities $m'_i$ of $\alpha_{I_{i_1}}^{e_{i_1}} \cdots \alpha_{I_{i_k}}^{e_{i_k}}$ satisfy $m'_i = \min\{m_i, 2\}$ and the second sum is over all remaining subsequences such that $m'_i = m_i$ if $m_i \leq 1$, $m'_i = 2$ if $m_i = 2$ and $m'_i = 1$ or $2$ if $m_i > 2$. The coefficients are defined by the formula: $c_{i_1 \cdots i_k} = \prod_r (m_{i_r} - 2)$, where the product is over all $r$ such that $m_{i_r} > 2$ and $m'_{i_r} = 1$. Note that there must be at least one such $r$.

Step 2. After Step 1 we can assume every $m_i \leq 2$. Now we replace every occurrence of $\alpha_I^{\pm 2}$ by $4\alpha_I^{\pm 1}$ and any occurrence of $\alpha_I \alpha_I^{-1}$ by 1.

Step 3. After Step 2 we can assume, in addition, that $I_i = I_j$ if and only if $i = j$. At this point any such monomial is determined by its corresponding oriented $UT$ graph. We describe the next two reductions in terms of these graphs and linear combinations of them. First we describe a preliminary modification. A *connected* $UT$ graph $G$ will be called *even* if every cycle has an even number of edges. In this case the vertices of $G$ can be divided into two classes: any two vertices in the same class are connected by a path with an even number of edges. If we choose one of these classes and remove each trivalent vertex in the class from the graph, replacing it by three univalent vertices, we obtain a new graph which is a union of $Y$ graphs and isolated edges. We carry over the orientations of $G$ to the new graph. Now define $G'$ to be the sum $G_1 + G_2$ of these new graphs, if $G$ is even, and 0 otherwise. We can now describe Step 3 by replacing every component $G$ with at least one univalent vertex, in an oriented $UT$ graph, by $G'$.

Step 4. After Step 3 we have a linear combination of oriented Ohtsuki graphs (with, perhaps, some additional isolated edges). This step will reduce it to a linear combination of Ohtsuki graphs with a preferred orientation at each trivalent vertex- for example, we can always prefer the orientation given by $i, j, k$ where $i < j < k$. We describe Step 4 as follows. Suppose $G$ is a component all of whose vertices are trivalent and $k$ of its vertices have the wrong orientation. Then we replace $G$ by $(-1)^k \tilde{G} + 2kG'$, where $\tilde{G} = G$ except that the orientations are now all the preferred ones. If $G$ is a $Y$ graph, then we replace $G$ by $\tilde{G}$.

Step 5. The final step is to eliminate any graph which has an isolated edge.

We leave as a straightforward, if, perhaps, lengthy exercise to show that these five steps can be achieved by using the relations given by equations 20-25 and 29. The important thing is to show that the result of this reduction depends only on the class of the monomial in $\mathbb{Q}[P_2^{\mathbb{SE}}(n)]/\mathcal{R}_n$. But this can be achieved by checking that, for each of the equations 20-23 and 29, reducing both sides of the equation gives the same result. □

**Corollary 3.5.** *$\mathcal{G}_m\mathcal{O}$ is a finite dimensional vector space and is nonzero only if $m$ is a multiple of 3.*

*Remark* 3.6. This partially answers question 1 of [Ga].

*Remark* 3.7. In [GO] we will introduce the notion of the 4 term relation and the weight system for a finite type invariant of $\mathbb{Z}HS$.

**3.4. From knots to 3-manifolds.** In this section we prove a vanishing theorem (theorem 7) for type $5m+1$ invariants of $\mathbb{Z}HS$. As a corollary, in proposition 3.9 we make some progress on question 2 of [Ga].

In order to simplify notation, for an $AS$ admissible link in $S^3$ we denote by $[L]$ the element $(S^3, L, \{+1, \dots, +1\})$ in $\mathcal{M}$ defined by equation (15).

**Theorem 7.** *If $L$ is an $AS$ admissible link in $S^3$ with a $(4m+1)$-component trivial sublink, then $[L] \in \mathcal{F}_{5m+2}\mathcal{M}$.*

*Proof.* If $L$ has $4m+1+r$ components, we proceed by downward induction on $r$. Obviously if $r > m$, there is nothing to prove. We record the following consequences of the defining relations in $\mathcal{M}$.

**Lemma 3.8.** *If $\tilde{L}$ is obtained from $L$ by changing a crossing of two bands, then $[\tilde{L}] - [L]$ is a linear combination of $[L_i]$, where each $L_i$ contains $L$ as a* proper *sublink.*

A *band* means a collection of parallel subarcs of components of $L$, and we assume, after choosing some oriention for the components of $L$ and a direction for each band, that each component of $L$ has an equal number of subarcs travelling in the positive direction in the two bands as in the negative direction. (Note that this condition is independent of the orientations chosen but does depend on the dirctions of the bands.) See figures 4 and 5.

*Proof.* It follows from figures 4, 5 and equations (15) and (19). □

Suppose we write $L = L_0 \cup L'$, where $L_0$ is the trivial sublink. If we perform a crossing change, as in Lemma 3.8 where at most one of the two bands contains arcs from $L_0$, then $\tilde{L}$ and each $L_i$ also contains $L_0$ as a sublink and so, by the inductive assumption, $[L] \in \mathcal{F}_{5m+2}\mathcal{M}$ if and only if $[\tilde{L}] \in \mathcal{F}_{5m+2}\mathcal{M}$. So, for example, we are free to change any component of $L'$ within its homotopy class in the complement of $L_0$. Furthermore, if $\pi = \pi_1(S^3 - L_0)$ then we can change a component $l$ by any element of $\pi_3$, by the argument in [Le1, pages 58-9], using band crossings in which one band consists of arcs from $L_0$ and the other consists of arcs from $l$.

As a consequence of these observations, we may assume that $L'$ consists of components $l_k$ so that each $l_k$ represents a product of commutators of degree two $\prod_{i<j}[x_i, x_j]^{e^k_{ij}}$. Now each commutator $[x_i, x_j]$ can be represented by a curve $\sigma_{ij}$ which intersects only two of the disjoint disks $D_j$ bounded by the components of $L_0$. See figure 15. We may therefore assume that each $l_k$ is a band sum of a number of copies of the $\sigma_{ij}$, slightly translated so that they are all disjoint. We now want to apply theorem 5. This tells us that $[L]$ is a linear combination of $[L_i]$, where $L_i$ coincides with $L$ except that each $l_k$ has been replaced by either one of the $\sigma_{ij}$ or a band sum of two of the $\sigma_{ij}$, or several of these. We can ignore any terms in which any $l_k$ has been replaced by more than one component, by the inductive assumption. In each of the $L_i$ that remain, each new $l_k$ intersects at most four of the $D_j$ and so $L_i'$ intersects at most $4r$ of the $D_j$. Therefore if $r \leq m$, there will be at least one of the $D_j$ not intersected by any of the component of $L$. But it then follows from equation (20) that $[L] = 0$. $\square$

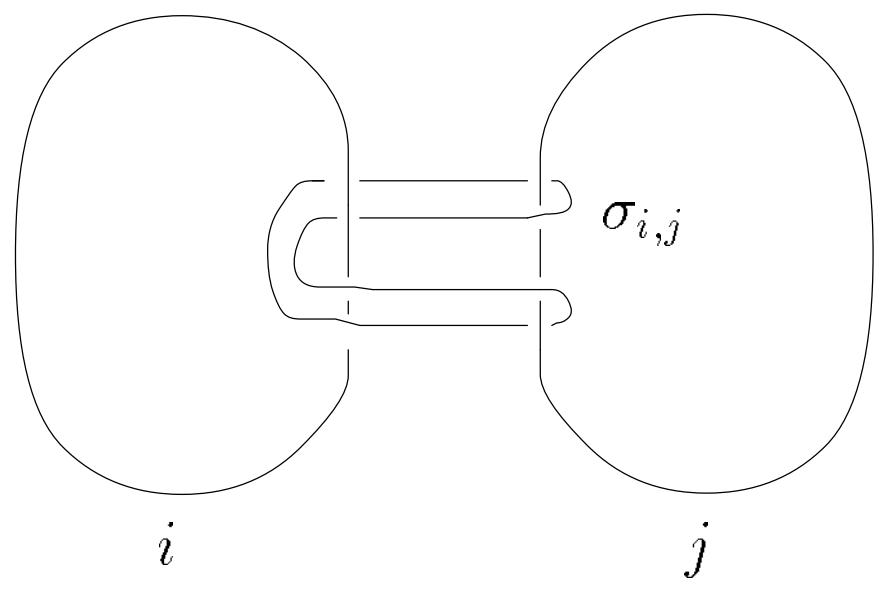

**Figure 15.** A curve $\sigma_{ij}$ that represents a commutator $[x_i, x_j]$ in the fundamental group of the complement of a trivial 2 component link.

As an application of Theorem 7, we study the map from knots to (linear combinations of) 3-manifolds defined by $K \to (S^3, K, +1) := (S^3)^{(K,+1)} - (S^3)$. Dually this map induces a map from 3-manifold invariants to knot invariants. In [Ga] we showed that this descends to a map

$$\mathcal{F}_n\mathcal{O} \longrightarrow \mathcal{F}_{n-1}\mathcal{V} \tag{30}$$

where $\mathcal{F}_n\mathcal{V}$ is the space of type $n$ knot invariants [B-N], [BL], [Va]. In [Ga] it is conjectured that the above map actually descends to a map:

$$\mathcal{F}_{3m}\mathcal{O} \longrightarrow \mathcal{F}_{2m}\mathcal{V} \tag{31}$$

**Proposition 3.9.** *The above map descends to a map:*

$$\mathcal{F}_{5m+1}\mathcal{O} \longrightarrow \mathcal{F}_{4m}\mathcal{V} \tag{32}$$

*Proof.* Let $\lambda \in \mathcal{F}_{5m}\mathcal{O}$ be a type $5m$ invariant of $\mathbb{Z}HS$, and let $\psi_\lambda$ be the associated knot invariant of equation (30). Let $K$ be an immersed knot in $S^3$ with $4m+1$ double points. Let $K \cup L_0$ denote the $AS$ admissible link in $S^3$ of $4m+2$ components obtained by replacing each double point with the left hand of figure 16. Note that $K \cup L_0$ contains an unlink $L_0$ of $4m+1$ components. Using the equality of figure 16 and the definition of the associated knot invariant $\psi_\lambda$ we obtain the following equality

$$\psi_\lambda(K) = \lambda([K \cup L_0]) \tag{33}$$

Using theorem 7 we have that $[K \cup L_0] \in \mathcal{F}_{5m+2}\mathcal{M}$ therefore $\lambda([K \cup L_0]) = 0$. This shows that $\psi_\lambda$ is a type $4m$ invariant of knots in $S^3$. □

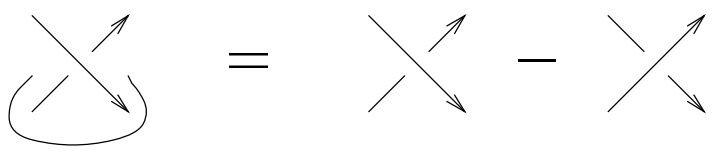

**Figure 16.** An equality in $\mathcal{M}$ related to crossing a change of a component of a link.

*Remark* 3.10. Proposition 3.9 is a $\frac{4m}{5m+1} \simeq \frac{4}{5}$ result, whereas question 2 of [Ga] asks for a $\frac{2}{3} < \frac{4}{5}$ result. In a recent preprint M. Greenwood and Xiao-Song Lin [GL] have shown a $\frac{n-2}{n} \simeq 1$ result, which is a weaker statement than proposition 3.9 if $n > 6$.

## 4. A philosophical comment

We feel that we owe a word on the appearance of trivalent graphs. As mentioned in the introduction, a motivation for the notion of finite type invariants of $\mathbb{Z}HS$ is Chern-Simons theory, exploited by Witten [Wi2]. Chern-Simons thoery is a *topological quantum field theory* with a topological Lagrangian containing a quadratic and a cubic term. The asymptotic expansion of the associated path integral (over the space of connections) as the coulping parameter goes to infinity can be approximated by a power series sum, each term of which is a finite sum over trivalent graphs (Feynman diagrams). This is the reason that trivalent graphs appear in Chern-Simons theory.

In the theory of finite type knot invariants, trivalent graphs appear either as triple point degenarations of knots, or as Feynman diagrams of an associated conformal field theory (governed by the $KZ$ equation), see [Dr], [B-N] and [Ko].

In the theory of finite type invariants of $\mathbb{Z}HS$ trivalent graphs appear because of the presence of the Kirby moves, an intrinsic 3-dimensional property of space.

Finally, in the notion of surgical equivalence of links trivalent graphs appear because of the association of figure 5 with triple commutators $[a, [b, c]]$ in fundamental groups.

The use of trivalent graphs, whether they come from the topology of 3-dimensional space, or the algebra (commutator groups, or cubic interaction terms in path integrals) is a unifying approach, and as such, it can be a source of inspiration, or confusion. We will let the reader decide which.

Department of Mathematics, Massachusetts Institute of Technology, Cambridge, MA 02139

*E-mail address*: stavros@math.mit.edu

Department of Mathematics, Brandeis University, Waltham, MA 02254-9110

*E-mail address*: levine@max.math.brandeis.edu